\documentclass[12pt]{amsart}
\usepackage{amscd,amssymb,times}
\newtheorem{thm}{Theorem}
\newtheorem{prop}[thm]{Proposition}
\newtheorem{lem}[thm]{Lemma}

\newtheorem{conj}[thm]{Conjecture}

\theoremstyle{remark}

\newtheorem{rem}[thm]{Remark}

\theoremstyle{definition}

\DeclareMathOperator{\Spc}{Spin^c}

\newcommand{\C}{\mathbb{ C}}

\newcommand{\Z}{\mathbb{ Z}}

\newcommand{\Diff}{\operatorname{Diff}}

\newcommand{\ra}{{\rightarrow}}

\begin{document}

\title{Dissolving four-manifolds and positive scalar curvature}
\author{B.~Hanke}
\address{Mathematisches Institut, Ludwig-Maximilians-Universit\"at
M\"unchen, Theresienstr.~39, 80333 M\"unchen, Germany}
\email{Bernhard.Hanke@mathematik.uni-muenchen.de}
\author{D.~Kotschick}
\email{dieter@member.ams.org}
\author{J.~Wehrheim}
\email{Jan.Wehrheim@mathematik.uni-muenchen.de}

\date{December 18, 2002; MSC 2000 classification: primary 57R57; 
secondary 53C25, 57R55}

\thanks{We are very grateful to R.~Gompf for valuable comments.
The first and the second authors are members of the
{\sl European Differential Geometry Endeavour} (EDGE),
Research Training Network HPRN-CT-2000-00101, supported by
The European Human Potential Programme}

\begin{abstract}
    We prove that many simply connected symplectic four-manifolds 
    dissolve after connected sum with only one copy of $S^{2}\times S^{2}$.

    For any finite group $G$ that acts freely on the three-sphere we 
    construct closed smooth four-manifolds with fundamental group $G$ 
    which do not admit metrics of positive scalar curvature, but whose 
    universal covers do admit such metrics. 
\end{abstract}

\maketitle

\section{Introduction}

It is a classical result of Wall~\cite{W} that any two simply connected 
four-manifolds with isomorphic intersection forms become diffeomorphic 
after taking the connected sum with sufficiently many copies of 
$S^{2}\times S^{2}$. It follows that any simply connected four-manifold 
is stably diffeomorphic to a connected sum of complex projective planes 
(with both orientations allowed) or to a connected sum of copies of 
$S^{2}\times S^{2}$ and of the $K3$ surface. In general, it is a 
very hard problem to determine the minimal number of copies of 
$S^{2}\times S^{2}$ required. Gauge theory shows that this number is 
usually positive, but, in spite of various attempts, it has not led to 
any further lower bounds, compare~\cite{D,FS}. This can be taken as 
evidence for the conjecture that one copy of $S^{2}\times S^{2}$ always 
suffices, which is only known in very few cases, for example for elliptic 
complex surfaces, see~\cite{Msurv}. 

One purpose of this paper is to prove that many simply connected 
symplectic four-manifolds constructed from algebraic surfaces by 
symplectic sums along submanifolds~\cite{Go1} become diffeomorphic to 
standard manifolds after taking the connected sum with only one copy 
of $S^{2}\times S^{2}$. We focus on the spin manifolds 
of nonnegative signature constructed by J.~Park~\cite{P}, but our 
argument applies to many other cases. There are several reasons for 
looking at these particular manifolds. Firstly, in the spin case there 
are gauge theoretic invariants~\cite{D,FS} which could, in theory, 
obstruct the kind of result we seek. Secondly, simply connected algebraic 
surfaces or symplectic manifolds with nonnegative signature are quite 
difficult to construct, and are considered to be more exotic than the 
ones of negative signature. 
Finally, we are interested in spin manifolds of zero signature because 
of an application to questions about the existence of positive scalar 
curvature metrics. 

In Section~\ref{s:pos} we give examples of both spin and non-spin 
four-manifolds with finite fundamental groups which do not admit metrics 
of positive scalar curvature, although their universal covers do admit 
such metrics. The possible fundamental groups are all finite groups 
which act freely on $S^{3}$, including all the finite cyclic groups. 
Previously, such an example with fundamental group of order two was given 
by LeBrun~\cite{LB}, whose result is analogous to one obtained by 
B\'erard~Bergery~\cite{BB} in high dimensions. The case of odd order 
fundamental groups is more interesting, because it disproves the 
following:
\begin{conj}[Rosenberg~\cite{Ro}, 1.2]\label{conj}
    Assume that $M^n$ is a connected closed manifold and $\pi_1(M)$ is
    finite of odd order. Then $M$ admits a metric of positive scalar
    curvature, if and only if its universal cover does.
\end{conj}
In fact, in dimensions $\geq 5$, Rosenberg~\cite{Ro} proved this 
conjecture for cyclic groups, and thus our examples with finite cyclic 
fundamental groups of odd orders do not have higher-dimensional analogs. 
Kwasik--Schultz~\cite{KS} have confirmed Rosenberg's conjecture for some 
classes of non-cyclic groups in dimensions $\geq 5$.

Since the beginning of Seiberg--Witten theory~\cite{Wi}, it has been 
known that there are many simply connected four-manifolds which do not 
admit metrics of positive scalar curvature, although the Gromov--Lawson 
conjecture~\cite{GL}, which is true in dimensions $\geq 5$, 
predicts that they should have such metrics. Our results in 
Section~\ref{s:pos} show that the situation is similar for 
Conjecture~\ref{conj}.

\section{Dissolving spin four-manifolds}\label{s:dissolve}

In this section we prove that there are many spin symplectic
four-manifolds which upon taking the connected sum with just one copy 
of $S^{2}\times S^{2}$ dissolve into connected sums of copies of 
$S^{2}\times S^{2}$ (in the zero signature case), or of copies of 
$S^{2}\times S^{2}$ and of the $K3$ surface (in the case of non-zero
signature). In fact, it will turn out that the requirement to 
dissolve after a single stabilization does not substantially 
influence the geography of Chern numbers.
See Remark~\ref{r:nonspin} for the non-spin case.

Our proof is based on the following result of Gompf~\cite{Go} 
(Lemma~4 and Corollary~5) elaborating on earlier work of 
Mandelbaum~\cite{Ma}. We only state the special case that we need.
\begin{prop}[Gompf~\cite{Go}]\label{p:diss}
    Let $M$ and $N$ be simply connected oriented $4$-manifolds
    containing the same embedded surface $F$ of genus $g\geq 1$ with 
    zero selfintersection. Assume that $F$ has simply connected 
    complement in $M$, and that $M$ is spin. Denote by $P$ the sum of $M$
    and $N$ along $F$. Then $P\# (S^{2}\times S^{2})$ is diffeomorphic
    to $M\# N\# 2g(S^{2}\times S^{2})$.
    \end{prop}
This implies in particular that the connected sum of any simply
connected spin elliptic surface with $S^{2}\times S^{2}$ is
diffeomorphic to a connected sum of copies of the $K3$ surface and of
$S^{2}\times S^{2}$, see~\cite{Ma}, or~\cite{Go} Corollary~8.

Here is our main result about spin manifolds of zero signature, which 
we will use for our application to questions about positive scalar 
curvature.
\begin{thm}\label{t:diss}
    There are infinitely many integers $l$ for which the
    following statements hold:
    \begin{enumerate}
    \item There are infinitely many symplectic manifolds $X_{i}$
    homeomorphic to the connected sum of $l$ copies of
    $S^{2}\times S^{2}$, which are pairwise non-diffeomorphic.
    \item The connected sum of each $X_{i}$ with $S^{2}\times S^{2}$
    is diffeomorphic to the connected sum of $l+1$ copies of
    $S^{2}\times S^{2}$.
    \end{enumerate}
    \end{thm}
\begin{proof}
    If one does not insist on the second property, the required
    examples have been constructed by J.~Park~\cite{P}. For certain $l$,
    Park first constructs a symplectic manifold $X$ homeomorphic to the
    connected sum of $l$ copies of $S^{2}\times S^{2}$. Then he shows
    using the knot surgery of Fintushel and Stern that one can change
    the smooth structure of $X$ to infinitely many different ones, all
    of which support symplectic structures.
    
    We will show that for infinitely many values of $l$ the first part
    of Park's construction yields a manifold $X$ which dissolves after
    connected sum with just one copy of $S^{2}\times S^{2}$. Then we
    construct infinitely many homeomorphic non-diffeomorphic examples
    $X_{i}$ from $X$ by varying one of the building blocks. Our
    construction allows us to check that all the $X_{i}$ do indeed
    dissolve upon connected sum with only one copy of
    $S^{2}\times S^{2}$.
    
    We now recall the construction of Park~\cite{P}. One begins with
    a simply connected spin algebraic surface $Y$ of positive
    signature containing a smooth holomorphic curve $F$ of genus $g$ of
    zero selfintersection, together with an embedded $2$-sphere $S$
    intersecting $F$ transversely in one point. Such examples are
    provided by Persson--Peters--Xiao~\cite{PPX}. The existence of $S$
    ensures that the complement of $F$ in $Y$ is simply connected.
    
    The second building block is a simply connected spin symplectic
    manifold $Z$ containing a symplectically embedded copy of $F$ also
    with zero selfintersection, and a symplectically embedded torus
    $T$ of zero selfintersection which is disjoint from $F$. Park
    exhibits concrete examples for each $g$.
    
    Let $X(k,n)$ be the manifold obtained by symplectically summing
    $k$ copies of $Y$ and one copy of $Z$ along $F$ and then
    symplectically summing the result with the simply connected spin
    elliptic surface $E(2n)$ without multiple fibers along $T$ and a
    fiber in $E(2n)$. The manifold $X(k,n)$ is spin and symplectic by
    construction. It is also simply connected because $F$ has simply
    connected complement in $Y$, as does the fiber in $E(2n)$. We
    shall think of $Z$ as the central component, with $k$ copies of $Y$
    and one copy of $E(2n)$ attached to $Z$, rather than being attached
    to each other. 
    
    By Novikov additivity the signature of $X(k,n)$ is the sum of the
    signatures of all the building blocks. Thus, if we choose
    \begin{equation}\label{n}
    n=\frac{1}{16}(k\sigma(Y)+\sigma(Z)) \ ,
    \end{equation}
    which is positive for large enough $k$ because $\sigma(Y)$ is
    positive, then $X(k,n)$ has zero signature. (Note that by
    Rochlin's theorem $\sigma(Y)$ and $\sigma(Z)$ are divisible by
    $16$.) By Freedman's classification~\cite{F}, $X(k,n)$ is then
    homeomorphic to a connected sum of $l$ copies of $S^{2}\times
    S^{2}$, where $l=\frac{1}{2}b_{2}(X(k,n))$.
    
    Now consider the connected sum $X(k,n)\# (S^{2}\times S^{2})$. By
    using Gompf's result, Proposition~\ref{p:diss} above, repeatedly
    we see that 
    $$
    X(k,n)\# (S^{2}\times S^{2}) \cong 
    kY\# Z\# E(2n)\# (2gk-k+2)(S^{2}\times S^{2}) \ . 
    $$
    As $E(2n)$ dissolves upon connected sum with $S^{2}\times S^{2}$,
    we obtain
    $$
    X(k,n)\# (S^{2}\times S^{2}) \cong 
    kY\# Z\# nK3\# (2gk-k+1+n)(S^{2}\times S^{2}) \ .
    $$
    Recall that $n$ grows linearly with $k$ according to~\eqref{n}.
    
    Wall~\cite{W} proved that any two simply connected $4$-manifolds 
    with isomorphic intersection forms become diffeomorphic after some 
    number of stabilizations with $S^{2}\times S^{2}$. Thus, if $k$ is 
    large enough, the connected sum of $Y$ and 
    $(2gk-k+1+n)(S^{2}\times S^{2})$ will be diffeomorphic to a connected 
    sum of copies of $S^{2}\times S^{2}$ and of copies of the $K3$ surface 
    with the non-complex orientation. We can then break up all the $k$ 
    copies of $Y$ inductively. Moreover, pairing the resulting copies of 
    the $K3$ surface with the non-complex orientation with copies of $K3$
    with the complex orientation, we obtain further copies of
    $S^{2}\times S^{2}$. The number of copies of $S^{2}\times S^{2}$
    that we split off grows with $k$, so by Wall's result we may assume 
    that there are enough of them to break up $Z$ as well. As $X(k,n)$ is 
    spin of zero signature, we finally see that for all large enough $k$ 
    and $n$ given by~\eqref{n}, the manifold $X(k,n)\# (S^{2}\times S^{2})$ 
    is diffeomorphic to a connected sum of copies of $S^{2}\times S^{2}$.
    
    It remains to show that there are infinitely many symplectic
    manifolds homeomorphic but non-diffeomorphic to $X(n,k)$, all of
    which dissolve upon connected sum with $S^{2}\times S^{2}$. For
    this we would like to replace the elliptic surface $E(2n)$ without 
    multiple fibers by one with multiple fibers obtained by logarithmic
    transformation. However, in this case the general fiber becomes 
    divisible in homology, in particular its complement is no longer simply
    connected. This could introduce fundamental group in our $4$-manifold, 
    and could obstruct the application of Proposition~\ref{p:diss}.
    
    To circumvent this problem we argue as follows. Think of the elliptic 
    surface $E(2n)$ as the fiber sum of $E(2n-2)$ with $E(2)=K3$. Instead 
    of attaching $E(2n)$ to $Z$ in the above construction, we attach $K3$ 
    to $Z$ and then attach $E(2n-2)$ to the $K3$ surface. However, these 
    attachments are not performed along parallel two-tori in $K3$. We 
    exploit the existence of two disjoint Gompf nuclei in $K3$, each of 
    which contains a torus of selfintersection zero and a transverse 
    $2$-sphere intersecting the torus once. Gompf proved in~\cite{Go1}, 
    Section~3, that the $K3$ surface has a symplectic structure which makes 
    the tori in two disjoint nuclei into symplectic submanifolds. Therefore, 
    the above construction, using summation along tori in different 
    nuclei, can be performed symplectically. Now, the boundary circle  
    of a normal disk in a tubular neighbourhood of such a torus is 
    null-homotopic along the transverse $2$-sphere inside the nucleus, 
    so that summing $K3$ to $Z$ (with copies of $Y$ attached) does not 
    introduce any fundamental group. Then, attaching the elliptic surfaces 
    with multiple fibers obtained by logarithmic transformations on 
    $E(2n-2)$ along a torus in a different nucleus of the $K3$ surface, 
    this torus will again have simply connected complement in its own 
    nucleus, making it irrelevant that the fiber of the elliptic 
    surface does not have simply connected complement. Thus 
    Proposition~\ref{p:diss} is applicable to all these summations.
    
    The logarithmic transformations on $E(2n-2)$ produce infinitely many
    distinct smooth structures on the topological manifold underlying
    $E(2n-2)$, which are detected by Seiberg--Witten invariants,
    cf.~\cite{Wi}. This difference in the Seiberg--Witten invariants
    survives the symplectic sum operation along a fiber, because of the
    gluing formulas due to Morgan--Mr\'owka--Szab\'o~\cite{MMS} and
    Morgan--Szab\'o--Taubes~\cite{MST}. Thus, we can produce infinitely
    many symplectic spin manifolds homeomorphic but non-diffeomorphic
    to $X(k,n)$. All these dissolve upon taking the connected sum with
    only one copy of $S^{2}\times S^{2}$ by the same argument as for
    $X(k,n)$.  
    \end{proof}


We can easily adapt the above argument to spin manifolds of 
nonzero signature. Starting from the examples constructed in the 
proof of Theorem~\ref{t:diss}, we can symplectically sum in extra 
copies of $Y$ along $F$ to make the signature positive, or we can use 
larger elliptic surfaces, with $n$ larger than the value given by~\eqref{n}, 
to make the signature negative. This will give manifolds which still 
dissolve after only one stabilization. 

To cover a large area of the geography more systematically, we 
proceed as follows. We use the coordinates $c_{1}^{2}=2e+3\sigma$ and 
$\chi =\frac{1}{4}(e+\sigma )$. The spin condition implies, via Rochlin's 
theorem, that 
\begin{equation}\label{spin}
c_{1}^{2}\equiv 8\chi\pmod{16} \ . 
\end{equation}
As spin symplectic manifolds are automatically minimal, they satisfy 
$c_{1}^{2}\geq 0$ by the work of Taubes~\cite{Taubes}. It is clear 
that in the simply connected case one must have $\chi > 0$. Thus, we 
try to cover lattice points in the first quadrant of the 
$(\chi , c_{1}^{2})$--plane subject to the congruence~\eqref{spin} with 
spin symplectic manifolds which dissolve after only one stabilization.

For a constant $c$ let $R_{c}$ denote the set of lattice points 
$(x,y)$ in the plane satisfying $y\equiv 8x\pmod{16}$, $x>0$, $y\geq 
0$, and 
\begin{alignat}{1}
    y &\leq 2x-16 \ , \\
    y &\leq 8x-c \ . \label{c}
    \end{alignat}

\begin{prop}\label{p:wedge}
There exists a constant $c$ such that all lattice points in $R_{c}$ 
are realized as the Chern invariants $(\chi,c_{1}^{2})$ of infinitely 
many pairwise nondiffeomorphic simply connected symplectic spin 
manifolds, all of which dissolve upon taking the connected sum with 
$S^{2}\times S^{2}$.       
\end{prop}
\begin{proof}
    The argument is modelled on the above proof of Theorem~\ref{t:diss} 
    and the proof of Proposition~1 in~\cite{P}.
    
    Let $H$ be a spin Horikawa surface with $c_{1}^{2}(H)=8$ and 
    $\chi(H)=7$. This is a symplectic genus two fibration over the 
    two-sphere; we fix a fiber $F$ of this fibration. The complement 
    of $F$ in $H$ contains the Milnor fiber of the $(2,3,7)$-singularity, 
    and one can find a torus $T$ of zero selfintersection with a transverse 
    $2$-sphere inside the Milnor fiber. We can deform the symplectic 
    structure of $H$ so that $T$ is a symplectic submanifold. Let $H(k,n)$ 
    be the manifold obtained by taking the fiber sum of $k$ copies of $H$ 
    along $F$ and then summing the resulting manifold with the elliptic 
    surface $E(2n)$ along $T$. This is a symplectic spin manifold, which is 
    simply connected because of the existence of the transverse $2$-spheres 
    for $F$ and $T$. We have $c_{1}^{2}(H(k,n))=16k-8$ and 
    $\chi(H(k,n))=8k+2n-1$. By varying $k\geq 1$ and $n\geq 1$, we can cover 
    all the lattice points in $R_{c}$ with odd $x$ using these manifolds.
    
    Now consider the connected sum $H(k,n)\# (S^{2}\times S^{2})$. By
    using Gompf's result, Proposition~\ref{p:diss}, repeatedly we see that 
    $$
    H(k,n)\# (S^{2}\times S^{2}) \cong 
    kH\# E(2n)\# (3k-1)(S^{2}\times S^{2}) \ . 
    $$
    As $E(2n)$ dissolves upon connected sum with $S^{2}\times S^{2}$,
    we obtain
    $$
    H(k,n)\# (S^{2}\times S^{2}) \cong 
    kH\# nK3\# (3k+n-2)(S^{2}\times S^{2}) \ .
    $$
    By the result of Wall~\cite{W} there is a $k_{0}$ such that 
    $H\# k_{0}(S^{2}\times S^{2})$ dissolves. Therefore, 
    $H(k,n)\# (S^{2}\times S^{2})$ dissolves as soon as $3k+n-2\geq 
    k_{0}$, which is equivalent to~\eqref{c} with $c=16k_{0}+32$.
    
    Finally, as $n\geq 1$, we can perform logarithmic transformations 
    on the elliptic building blocks to achieve infinitely many distinct 
    smooth structures on the topological manifold underlying $H(k,n)$. 
    All these smooth structures admit symplectic structures and 
    still dissolve after only one stabilization. This completes the 
    proof for odd values of $x$.
    
    To cover the lattice points with even $x$ in $R_{c}$, we use 
    $H'(k,n)$, obtained from the above $H(k,n)$ by summing in an 
    additional copy of $H$, summed along $T$, not along $F$. The 
    resulting manifolds are again simply connected because $H$ 
    contains a transverse $2$-sphere for $T$ in the complement of $F$. 
    This construction covers all lattice points with $y>0$. For $y=0$ 
    we can just use the spin elliptic surfaces $E(2n)$ themselves. The 
    rest of the argument is as for the case of odd $x$.
    \end{proof}

    This leads to the following geography result:
\begin{thm}\label{t:pos}
    There is a line of slope $>8$ in the $(x,y)$-plane such that 
    every lattice point in the first quadrant which is below this 
    line and satisfies $y\equiv 8x\pmod{16}$ is realized by the Chern 
    invariants $(\chi,c_{1}^{2})$ of infinitely many pairwise
    nondiffeomorphic simply connected symplectic spin manifolds, all 
    of which dissolve upon taking the connected sum with $S^{2}\times S^{2}$.
    \end{thm}
    \begin{proof}
        Let $Y$ and $Z$ be the building blocks from the proof of 
        Theorem~\ref{t:diss}. We sum $k$ copies of $Y$ to $Z$ along the 
        surface $F$ of genus $g$. Then we sum the $H(l,n)$ and $H'(l,n)$ 
	from the region $R_{c}$ in Proposition~\ref{p:wedge} to the
	resulting manifolds by summation along the torus $T$ in $Z$ and 
	in the elliptic piece of $H(l,n)$ or $H'(l,n)$. In all these 
	summations the complement of the surface along which the summation 
	is performed is simply connected in at least one of the summands, 
	so that the resulting manifolds are simply connected.
        
        By the proofs of Theorem~\ref{t:diss} and of Proposition~\ref{p:wedge}, 
	these manifolds have all the desired properties as soon as $3l+n$ is 
	large enough. It is easy to see that varying $k$ and letting $(l,n)$ 
	range over the parameters of the $H(l,n)$ or $H'(l,n)$ in $R_{c}$, 
	the Chern invariants cover all the lattice points in the claimed area, 
	because $Y$ has positive signature. (We may have to increase the 
	constant $c$ from Proposition~\ref{p:wedge} in order to ensure that 
	$3l+n$ is always large enough.)
        \end{proof}
\begin{rem}
    Theorem~\ref{t:pos} should be compared to the main theorem of 
    Park~\cite{P}, who proved a version of it without the requirement that 
    the manifolds in question dissolve after only one stabilization. 
    Park argues that by using many copies of $Y$, one can push the slope up 
    to approximately the slope $\frac{c_{1}^{2}(Y)}{\chi(Y)}$ of $Y$, 
    and~\cite{PPX} provides a construction for $Y$ with slope $> 8.76$. 
    However, the summation of $Y$ to itself, or to $Z$, is performed along 
    a surface of (unknown) genus $g>1$, so that the Chern numbers are not 
    additive. Instead, asymptotically for large $k$, the best slope one can 
    obtain is approximately
    $$
    \frac{c_{1}^{2}(Y)+8(g-1)}{\chi(Y)+g-1} > 8 \ ,
    $$
    both in Park's result and in ours. This is smaller than 
    $\frac{c_{1}^{2}(Y)}{\chi(Y)}$.
    
    Thus, considering only manifolds which dissolve after a single 
    stabilization does not alter the geography in any essential way.
    \end{rem}
    

\begin{rem}
    For the spin manifolds in Theorems~\ref{t:diss} and~\ref{t:pos}, 
    the $2$-torsion instanton invariants of~\cite{D,FS} are defined, but 
    vanish because the manifolds dissolve after only one stabilization. 
    Therefore~\cite{FS}, their Donaldson polynomials are all even. Under 
    certain technical hypotheses, Fintushel--Stern~\cite{FS} proved the 
    evenness of Donaldson polynomials for manifolds which do not necessarily 
    dissolve after the first stabilization.
    \end{rem}
    
\begin{rem}\label{r:nonspin}
    Most of our examples also have the property that they dissolve upon
    connected sum with a single copy of $\C P^{2}$. This means that 
    they are almost completely decomposable in the sense of 
    Mandelbaum~\cite{Msurv}. These are the first examples of irreducible 
    four-manifolds of non-negative signature which are almost completely 
    decomposable.

    One can use the arguments in the proofs of Theorems~\ref{t:diss} 
    and~\ref{t:pos} to exhibit non-spin almost completely decomposable minimal 
    symplectic manifolds of positive and of zero signature by using non-spin 
    elliptic surfaces instead of spin ones, or by constructing similar 
    irreducible manifolds starting from non-spin building blocks, rather 
    than the spin ones of Persson-Peters-Xiao~\cite{PPX}. The geography 
    statements one obtains for irreducible almost completely decomposable 
    four-manifolds in the non-spin case are rather stronger than 
    Theorem~\ref{t:pos}. 
    \end{rem}

\section{Positive scalar curvature and finite coverings}\label{s:pos}

Let $M$ be a smooth Riemannian manifold of positive scalar curvature. 
Clearly, if $\bar M \rightarrow M$ is an unramified covering of
$M$, then the pulled back metric on $\bar M$ also has positive
scalar curvature. On the other hand, if $M$ is a smooth manifold
and a finite cover $\bar M$ of $M$ admits a Riemannian metric of
positive scalar curvature, then it is not usually true that $M$ also
admits such a metric. One might try to average the metric on $\bar M$ 
and consider the induced metric on $M$, but this approach has turned
out to be too naive. Using an index-theoretic obstruction with
values in the $K$-theory of a certain $C^*$-algebra associated to
the fundamental group of the manifold under consideration,
Rosenberg~\cite{Ro1} exhibited a triple cover $\bar M \ra M$ of a
closed five-dimensional manifold for which $\bar M$ admits a metric
of positive scalar curvature, but $M$ does not. In this example $M$
has infinite fundamental group. Rosenberg pointed out in the same
paper that the situation is very different for manifolds with finite 
fundamental groups. 
Examples constructed by B\'erard Bergery~\cite{BB} 
show that there are 
high-dimensional closed smooth manifolds with fundamental group
$\Z/2$ which do not admit metrics of positive scalar curvature, 
although their universal covers do. This led Rosenberg~\cite{Ro} to 
formulate Conjecture~\ref{conj} for odd order finite fundamental groups.

We now disprove this conjecture in dimension $4$. More generally, we 
show the following:
\begin{thm}\label{counter}
    For any nontrivial finite group $G$ which acts freely on $S^{3}$,
    there are closed smooth four-manifolds $M$ with fundamental group 
    $G$ which do not admit metrics of positive scalar curvature, but 
    whose universal covers do admit such metrics. The manifolds $M$ 
    can be taken to be either spin or non-spin.
\end{thm}
Every cyclic group acts freely on $S^{3}$ and among odd order groups 
there are no others by Hopf's theorem~\cite{Hopf}. Thus, 
Theorem~\ref{counter} disproves Rosenberg's conjecture exactly for 
all finite cyclic groups of odd order. Moreover, there are both spin 
and non-spin counterexamples. The non-spin case of Theorem~\ref{counter}
with $G=\Z/2$, which is not relevant to Rosenberg's conjecture, was
proved by LeBrun~\cite{LB}. His argument works only for non-spin
manifolds, whereas the one of B\'erard Bergery~\cite{BB} in high 
dimensions works only for spin manifolds.

\begin{proof}[Proof of Theorem~\ref{counter}]
Let $G$ be a finite group of order $d>1$ acting freely on $S^{3}$, 
and let $L$ be the quotient $S^{3}/G$. On the product $L\times S^{1}$ 
one can perform surgery to kill the fundamental group of the second 
factor in such a way that the resulting $4$-manifold $N$ is spin. It
is obviously a rational homology sphere.

The following is well-known, compare Proposition~1-3 of~\cite{Ue}, or~\cite{GS}.
\begin{lem}\label{diff}
    The universal cover $\tilde{N}$ of $N$ is diffeomorphic to
    the connected sum of $d-1$ copies of $S^2 \times S^2$.
\end{lem}

We first exhibit the claimed non-spin examples.
Let $X=E(2n+1)$ be a non-spin simply connected elliptic surface 
without multiple fibers, for some $n\geq 1$. Then the positive part 
of its second Betti number $b_2^+(X) > 1$. Therefore~\cite{Wi}, the 
Seiberg--Witten invariant of $X$ is well-defined and non-zero for the 
$\Spc$-structure induced by the K\"ahler structure of $X$. We set 
$M = N \# X$. This is a smooth closed $4$-manifold with fundamental 
group $G$ which is not spin. The following is an immediate consequence 
of the gluing result of~\cite{KMT}.
\begin{prop}[Kotschick-Morgan-Taubes~\cite{KMT}]\label{p:KMT}
    For every $\Spc$-structure on $X$ there is one on $M$ with the 
    same Seiberg--Witten invariant. In particular, $M$ does not 
    admit a metric of positive scalar curvature.
\end{prop}

By Lemma~\ref{diff} above, the universal covering of $M$ is 
diffeomorphic to $(d-1)(S^2 \times S^2) \# d \ X$. Now, the 
connected sum of any simply connected elliptic surface with
$S^2 \times S^2$ dissolves by an application of Proposition~\ref{p:diss},
see~\cite{Go,Msurv}. As $X$ is not spin, we see that the universal 
cover of $M$ is diffeomorphic to a connected sum of copies of 
$\C P^2$ and of $\overline{\C P^2}$.

By~\cite{GL,SY} the class of manifolds of dimension at least three
admitting  metrics of positive scalar curvature is closed under 
forming connected sums. Thus connected sums of copies of $\C P^2$ and 
of $\overline{\C P^2}$ admit such metrics.

In order to obtain spin examples for $M$, we need to replace $X$
in the above construction with a spin manifold with non-trivial
Seiberg--Witten invariants. Then the connected sum $M=N\# X$ and its
universal covering will also be spin. The signature of $\tilde M$ will
be $d$ times the signature of $M$, which is the same as the signature
of $X$. Therefore, the Lichnerowicz vanishing theorem forces us to
choose $X$ so that it has zero signature. We take for $X$ one of the
symplectic spin manifolds with zero signature constructed in the proof
of Theorem~\ref{t:diss}. Because this becomes diffeomorphic to a
connected sum of copies of $S^{2}\times S^{2}$ after only one
stabilization, the universal covering of $M$ is also a
connected sum of copies of $S^{2}\times S^{2}$, and thus admits
metrics of positive scalar curvature.

This completes the proof of Theorem~\ref{counter}.
\end{proof}

All our counterexamples, $M$, to Rosenberg's conjecture have the
additional property that they are homeomorphic to manifolds $M'$ which
do satisfy the conjecture, in that both $M'$ itself and its universal
covering $\tilde{M'}$ admit a metric of positive scalar curvature.
In fact, the universal covers $\tilde{M}$ and $\tilde{M'}$ are
diffeomorphic, and have standard differentiable structures
(connected sums of $\C P^{2}$ and $\overline{\C P^{2}}$, or of copies
of $S^{2}\times S^{2}$). The group $G$ acts freely on these standard
manifolds in two essentially different ways: there is the standard
action by isometries of a positive scalar curvature metric, with
quotient $M'$, and there are exotic actions, which do not fix
any positive scalar curvature metric, with exotic quotients $M$.
That gauge theory detects exotic group actions, though not the
application to positive scalar curvature metrics, was observed 
before, for example by Ue~\cite{Ue}.

We can use these group actions to say something about the space of
metrics of positive scalar curvature on certain manifolds.
\begin{thm}
    Let $G$ be any nontrivial finite group which acts freely on $S^{3}$. 
    Then $G$ acts freely by diffeomorphisms in infinitely many ways on
    infinitely many spin and infinitely many non-spin four-manifolds
    $M$ of positive scalar curvature, such that in each case all the 
    actions are conjugate by homeomorphisms but are not conjugate by 
    diffeomorphisms of $M$. These actions give rise to infinitely many 
    actions of $G$ without fixed points on the space of positive scalar 
    curvature metrics on $M$, which are not conjugate in $\Diff(M)$.
    \end{thm}
    In each case there is also one action which has a fixed point in
    the space of metrics of positive scalar curvature, giving rise to
    a standard quotient. The exotic actions have quotients without
    positive scalar curvature, and so have no fixed points on the
    space of metrics of positive scalar curvature. Note that if $G$
    is of prime order, the exotic actions on the space of positive
    scalar curvature metrics are free.
\begin{proof}
    In the spin case, the infinitely many examples are connected sums
    of different numbers of copies of $S^{2}\times S^{2}$. The group
    $G$ has infinitely many non-conjugate actions on a fixed such
    manifold giving rise to quotients of the form $X_{i}\# N$, with
    $X_{i}$ as in Theorem~\ref{t:diss}.
    
    For the non-spin case the $G$-actions are on connected sums of
    copies of $\C P^{2}$ and of $\overline{\C P^{2}}$, with quotients
    of the form $X_{j}\# N$, with $X_{j}$ infinitely many homeomorphic
    but pairwise non-diffeomorphic elliptic surfaces (obtained from each
    other by logarithmic transformation).
    
    Proposition~\ref{p:KMT} shows that the distinct differentiable
    structures on $X_{i}$ or $X_{j}$, which are detected by
    Seiberg--Witten invariants, remain distinct after connected sum
    with $N$.
    \end{proof}

\begin{rem}
    Beyond the non-spin case with $G$ of order two of the above
Theorem~\ref{counter}, LeBrun's paper~\cite{LB} contains a result
(Theorem 2) on manifolds with infinite fundamental groups, where he
considers finite coverings of high degree. The claimed coverings do
not exist, as they would violate the multiplicativity of the Euler
characteristic in finite coverings. However, the argument can be
salvaged by considering a correctly chosen sequence of coverings of
arbitrarily large degree. These coverings will always have even 
degrees.
\end{rem}

\end{document}